\documentclass[a4paper,11pt,oneside]{amsart}

\RequirePackage{amsmath}
\RequirePackage{bm}
\RequirePackage{amssymb}
\RequirePackage{upref}
\RequirePackage{amsthm}
\RequirePackage{enumerate}
\RequirePackage{pb-diagram}
\RequirePackage{amsfonts}
\RequirePackage[mathscr]{eucal}
\RequirePackage{verbatim}
\RequirePackage{xr}
\RequirePackage{graphicx}
\usepackage{calc}
\usepackage{xspace}

\newcommand{\al}{\alpha}
\newcommand{\bet}{\beta}
\newcommand{\ga}{\gamma}
\newcommand{\de}{\delta }
\newcommand{\e}{\epsilon}

\newcommand{\f}{\varphi}
\newcommand{\h}{\eta}

\newcommand{\ka}{\kappa}
\newcommand{\lam}{\lambda}
\newcommand{\m}{\mu}
\newcommand{\n}{\nu}

\newcommand{\s}{\sigma}
\newcommand{\x}{\xi}

\newcommand{\D}{\varDelta}

\newcommand{\Om}{\varOmega}

\newcommand{\di}[1]{#1\nobreakdash-\hspace{0pt}dimensional}
\newcommand{\nbdd}{\nobreakdash--}
\newcommand{\nbd}{\nobreakdash-\hspace{0pt}}

\newcommand{\fu}[3]{#1\hspace{0pt}_{|_{#2_#3}}}
\newcommand{\fv}[2]{#1\hspace{0pt}_{|_{#2}}}

\newcommand{\so}{{\mc S_0}}

\newcommand{\const}{\tup{const}}
\newcommand{\ndash}{\nobreakdash--}

\newcommand{\msp[1]}[1]{\mspace{#1mu}}
\newcommand{\low}[1]{{\hbox{}_{#1}}}

\newcommand{\R}[1][n+1]{{\protect\mathbb R}^{#1}}

\newcommand{\N}{{\protect\mathbb N}}

\newcommand{\eR}{\stackrel{\lower1ex \hbox{\rule{6.5pt}{0.5pt}}}{\msp[3]\R[]}}
\newcommand{\eN}{\stackrel{\lower1ex \hbox{\rule{6.5pt}{0.5pt}}}{\msp[1]\N}}
\newcommand{\eO}{\stackrel{\lower1ex
\hbox{\rule{6pt}{0.5pt}}}{\msc O}}

\DeclareMathOperator{\graph}{graph}

\newcommand\pa{\partial}
\newcommand\pde[2]{\frac {\partial#1}{\partial#2}}
 
\newcommand\sql[1][u]{\sqrt{1-|D#1|^2}}
\newcommand{\un}{\infty}
\newcommand{\A}{\forall}

\newcommand{\set}[2]{\{\,#1\colon #2\,\}}
\newcommand{\uu}{\cup}

\newcommand{\uuu}{\bigcup}

\newcommand{\uud}{ \stackrel{\lower 1ex \hbox {.}}{\uu}}
\newcommand{\uuud}[1]{ \stackrel{\lower 1ex \hbox {.}}{\uuu_{#1}}}
\newcommand\su{\subset}

\newcommand{\sminus}[1][28]{\raise 0.#1ex\hbox{$\scriptstyle\setminus$}}

\newcommand{\abs}[1]{\lvert#1\rvert}

\newcommand{\norm}[1]{\lVert#1\rVert}

\newcommand{\nnorm}[1]{| \mspace{-2mu} |\mspace{-2mu}|#1| \mspace{-2mu}
|\mspace{-2mu}|}
\newcommand{\spd}[2]{\protect\langle #1,#2\protect\rangle}

\newcommand\ch[3]{\varGamma_{#1#2}^#3}
\newcommand\cha[3]{{\bar\varGamma}_{#1#2}^#3}

\newcommand{\riem}[4]{R_{#1#2#3#4}}
\newcommand{\riema}[4]{{\bar R}_{#1#2#3#4}}

\newcommand{\tbf}{\textbf}
\newcommand{\tit}{\textit}

\newcommand{\tup}{\textup}

\newcommand{\mc}{\protect\mathcal}
\newcommand{\msc}{\protect\mathscr}

\providecommand{\bysame}{\makebox[3em]{\hrulefill}\thinspace}

\newcommand{\ci}{\cite}
\newcommand{\bib}{\bibitem}

\newcommand{\bt}{\begin{thm}}
\newcommand{\bl}{\begin{lem}}
\newcommand{\bc}{\begin{cor}}
\newcommand{\bd}{\begin{definition}}
\newcommand{\bpp}{\begin{prop}}
\newcommand{\br}{\begin{rem}}
\newcommand{\bn}{\begin{note}}
\newcommand{\be}{\begin{ex}}
\newcommand{\bes}{\begin{exs}}
\newcommand{\bb}{\begin{example}}
\newcommand{\bbs}{\begin{examples}}
\newcommand{\ba}{\begin{axiom}}

\newcommand{\et}{\end{thm}}
\newcommand{\el}{\end{lem}}
\newcommand{\ec}{\end{cor}}
\newcommand{\ed}{\end{definition}}
\newcommand{\epp}{\end{prop}}
\newcommand{\er}{\end{rem}}
\newcommand{\en}{\end{note}}
\newcommand{\ee}{\end{ex}}
\newcommand{\ees}{\end{exs}}
\newcommand{\eb}{\end{example}}
\newcommand{\ebs}{\end{examples}}
\newcommand{\ea}{\end{axiom}}

\newcommand{\bp}{\begin{proof}}
\newcommand{\ep}{\end{proof}}
\newcommand{\eps}{\renewcommand{\qed}{}\end{proof}}

\newcommand{\bal}{\begin{align}}

\newcommand{\bi}[1][1.]{\begin{enumerate}[\upshape #1]}
\newcommand{\bia}[1][(1)]{\begin{enumerate}[\upshape #1]}
\newcommand{\bin}[1][1]{\begin{enumerate}[\upshape\bfseries #1]}
\newcommand{\bir}[1][(i)]{\begin{enumerate}[\upshape #1]}
\newcommand{\bic}[1][(i)]{\begin{enumerate}[\upshape\hspace{2\cma}#1]}
\newcommand{\bis}[2][1.]{\begin{enumerate}[\upshape\hspace{#2\parindent}#1]}
\newcommand{\ei}{\end{enumerate}}

\newcommand\ndots{\raise 0.47ex \hbox {,}\hskip0.06em\cdots %
     \raise 0.47ex \hbox {,}\hskip0.06em} 

\newcommand{\q}{\quad}
\newcommand{\qq}{\qquad}

\newcommand\nd{\noindent}

\newskip\Csmallskipamount                                                
\Csmallskipamount=\smallskipamount
\newskip\Cmedskipamount
\Cmedskipamount=\medskipamount
\newskip\Cbigskipamount
\Cbigskipamount=\bigskipamount

\newcommand\cvs{\vspace\Csmallskipamount}   
\newcommand\cvm{\vspace\Cmedskipamount}

\newskip\csa
\csa=\smallskipamount

\newskip\cma
\cma=\medskipamount

\newskip\cba
\cba=\bigskipamount

\newdimen\spt
\newcommand\citem{\cvs\advance\itemno by
1{(\romannumeral\the\itemno})\hskip3pt}
\newcommand{\bitem}{\cvm\nd\advance\itemno by
1{\bf\the\itemno}\hspace{\cma}}

\newcount\itemno
\newcommand{\las}[1]{\label{S:#1}}

\newcommand{\lae}[1]{\label{E:#1}}
\newcommand{\lat}[1]{\label{T:#1}}
\newcommand{\lal}[1]{\label{L:#1}}

\newcommand{\lap}[1]{\label{P:#1}}
\newcommand{\lar}[1]{\label{R:#1}}

\newcommand{\rs}[1]{Section~\ref{S:#1}}

\newcommand{\rt}[1]{Theorem~\ref{T:#1}}
\newcommand{\rl}[1]{Lemma~\ref{L:#1}}

\newcommand{\rp}[1]{Proposition~\ref{P:#1}}
\newcommand{\rr}[1]{Remark~\ref{R:#1}}

\newcommand{\re}[1]{\eqref{E:#1}}

\newskip\thmskip
\thmskip=\parindent

\newskip\hsk
\newenvironment{hinw}{\labelsep=0pt\begin{list}{}{\labelsep=0pt\itemindent=0pt\labelwidth=0pt\leftmargin=\parindent\rightmargin=0pt\partopsep=\cba}%
\item\it\nopagebreak\nopagebreak}%
{\end{list}}

\newcommand\bh{\begin{hinw}}
\newcommand{\eh}{\end{hinw}}

\newtheoremstyle{normal}
  {\cba}
  {\cba}
  {}
  {\thmskip}
  {\bfseries}
  {.}
  {\hsk}
  {}

\newtheoremstyle{abschnitt}
  {\cba}
  {\cba}
  {}
  {\thmskip}
  {\bfseries}
  {.}
  {\hsk}
  {}

\newtheoremstyle{italic}
  {\cba}
  {\cba}
  {\itshape}
  {\thmskip}
  {\bfseries}
  {.}
  {\hsk}
  {}

\newtheoremstyle{aufgaben}
  {\cba}
  {\cba}
  {}
  {}
  {\normalsize\bfseries}
  {.}
  {\hsk}
  {}

\newtheoremstyle{break}
  {\cba}
  {\cba}
  {\itshape}
  {}
  {\bfseries}
  {.}
  {\newline}
  {}

\swapnumbers
\theoremstyle{italic}
\newtheorem{thm}[subsection]{Theorem}
\newtheorem{lem}[subsection]{Lemma}
\newtheorem{prop}[subsection]{Proposition}
\newtheorem{cor}[subsection]{Corollary}

\theoremstyle{normal}
\newtheorem{rem}[subsection]{Remark}
\newtheorem{definition}[subsection]{Definition}
\newtheorem{example}[subsection]{Example}
\newtheorem{examples}[subsection]{Examples}
\newtheorem{ex}[subsection]{Exercise}
\newtheorem{note}[subsection]{}
\newtheorem{axiom}[subsection]{Axiom}

\theoremstyle{aufgaben}
\newtheorem{exs}[subsection]{Exercises}

\swapnumbers

\numberwithin{equation}{section}
\numberwithin{figure}{section}

\newenvironment{textequation}[1][0.8]
{\begin{equation}
\begin{aligned}
\begin{minipage}{#1\linewidth}}
{\end{minipage}
\end{aligned}
\end{equation}
\ignorespacesafterend}

\newcommand{\btext}{\begin{textequation}}
\newcommand{\etext}{\end{textequation}}

\RequirePackage{upref}
\RequirePackage{amsthm}
\RequirePackage{enumerate}
\usepackage{xr-hyper}
\usepackage{hyperref}

\begin{document}
\title[Hypersurfaces of prescribed mean curvature]{Hypersurfaces of prescribed
mean curvature in Lorentzian manifolds}

\author{Claus Gerhardt}
\address{Ruprecht-Karls-Universit\"at, Institut f\"ur Angewandte Mathematik,
Im Neuenheimer Feld 294, 69120 Heidelberg, Germany}
\email{gerhardt@math.uni-heidelberg.de}

%
\subjclass{}
\keywords{Prescribed mean curvature, Lorentz manifold}
\date{April 12, 1999}
%


\begin{abstract} We give a new existence proof for closed hypersurfaces of
prescribed mean curvature in Lorentzian manifolds.
\end{abstract}
\maketitle

\tableofcontents
\setcounter{section}{-1}
\section{Introduction} 

Hypersurfaces of prescribed mean curvature especially those with constant
mean curvature play an important role in general relativity. In \ci{cg83} the
existence of closed hypersurfaces of prescribed mean curvature in a globally
hyperbolic Lorentz manifold with a compact Cauchy hypersurface was proved
provided there were barriers. The proof consisted of two parts, the a priori
estimates for the gradient and the application of a fixed point theorem. That
latter part of the proof was rather complicated, and certainly nobody would have
qualified it as elegant.

Ecker and Huisken, therefore, gave another existence proof using an
evolutionary approach, but they had to assume that the time-like convergence
condition is satisfied, and, even more important, that the prescribed mean
curvature satisfies a structural monotonicity condition, cf. \ci{eh91}. These are
serious restrictions which had to be assumed because the authors relied on the
gradient estimate of Bartnik \ci{rb84}, who had proved another a priori estimate
in the elliptic case.

We shall show in the following that the evolutionary method can be used in the
existence proof without any unnecessary restrictions on the curvature of the
ambient space or the right-hand side. The only difference in the
assumptions---relative to our former paper--- is that the right-hand side is
now supposed to be of class $C^1$, while bounded is actually sufficient. But this
drawback can easily be overcome by approximation.

This paper is organized as follows: In \rs{1} we introduce the notations and
definitions we rely on.

In \rs 2 we look at the curvature flow associated with our problem, and the
corresponding evolution equations for the basic geometric quantities of the
flow hypersurfaces.

In \rs 3  lower order estimates for the evolution problem are proved, while a
priori estimates in the $C^2$\nbd{norm} are derived in \rs 4.

Finally, in \rs 5, we demonstrate that the evolutionary solution converges to a
stationary solution.

\section{Notations and definitions}\las{1}

The main objective of this section is to state the equations of Gau{\ss}, Codazzi,
and Weingarten for hypersurfaces $M$ in a  \di{(n+1)} Lorentzian  space $N$.
Geometric quantities in $N$ will be
denoted by
$(\bar g_{\al\bet}),(\riema \al\bet\ga\de)$, etc., and those in $M$ by $(g_{ij}), (\riem
ijkl)$, etc. Greek indices range from $0$ to $n$ and Latin from $1$ to $n$; the
summation convention is always used. Generic coordinate systems in $N$ resp.
$M$ will be denoted by $(x^\al)$ resp. $(\x^i)$. Covariant differentiation will
simply be indicated by indices, only in case of possible ambiguity they will be
preceded by a semicolon, i.e. for a function $u$ in $N$, $(u_\al)$ will be the
gradient and
$(u_{\al\bet})$ the Hessian, but e.g., the covariant derivative of the curvature
tensor will be abbreviated by $\riema \al\bet\ga{\de;\e}$. We also point out that
\begin{equation}
\riema \al\bet\ga{\de;i}=\riema \al\bet\ga{\de;\e}x_i^\e
\end{equation}
with obvious generalizations to other quantities.

Let $M$ be a \tit{space-like} hypersurface, i.e. the induced metric is Riemannian,
with a differentiable normal $\n$ that is time-like.

In local coordinates, $(x^\al)$ and $(\x^i)$, the geometric quantities of the
space-like hypersurface $M$ are connected through the following equations
\begin{equation}\lae{1.3}
x_{ij}^\al= h_{ij}\n^\al
\end{equation}
the so-called \tit{Gau{\ss} formula}. Here, and also in the sequel, a covariant
derivative is always a \tit{full} tensor, i.e.

\begin{equation}
x_{ij}^\al=x_{,ij}^\al-\ch ijk x_k^\al+\cha \bet\ga\al x_i^\bet x_j^\ga.
\end{equation}
The comma indicates ordinary partial derivatives.

In this implicit definition the \tit{second fundamental form} $(h_{ij})$ is taken
with respect to $\n$.

The second equation is the \tit{Weingarten equation}
\begin{equation}
\n_i^\al=h_i^k x_k^\al,
\end{equation}
where we remember that $\n_i^\al$ is a full tensor.

Finally, we have the \tit{Codazzi equation}
\begin{equation}
h_{ij;k}-h_{ik;j}=\riema\al\bet\ga\de\n^\al x_i^\bet x_j^\ga x_k^\de
\end{equation}
and the \tit{Gau{\ss} equation}
\begin{equation}
\riem ijkl=- \{h_{ik}h_{jl}-h_{il}h_{jk}\} + \riema \al\bet\ga\de x_i^\al x_j^\bet x_k^\ga
x_l^\de.
\end{equation}

Now, let us assume that $N$ is a globally hyperbolic Lorentzian manifold with a
\tit{compact} Cauchy surface. $N$ is then a topological product $\R[]\times \mc
S_0$, where $\mc S_0$ is a compact Riemannian manifold, and there exists a
Gaussian coordinate system
$(x^\al)$, such that the metric in $N$ has the form
\begin{equation}\lae{1.7}
d\bar s_N^2=e^{2\psi}\{-{dx^0}^2+\s_{ij}(x^0,x)dx^idx^j\},
\end{equation}
where $\s_{ij}$ is a Riemannian metric, $\psi$ a function on $N$, and $x$ an
abbreviation for the space-like components $(x^i)$, see \ci{GR},
\ci[p.~212]{HE}, \ci[p.~252]{GRH}, and \ci[Section~6]{cg83}.
 We also assume that
the coordinate system is \tit{future oriented}, i.e. the time coordinate $x^0$
increases on future directed curves. Hence, the \tit{contravariant} time-like
vector$(\x^\al)=(1,0,\dotsc,0)$ is future directed as is its \tit{covariant} version
$(\x_\al)=e^{2\psi}(-1,0,\dotsc,0)$.

Let $M=\graph \fv u\so$ be a space-like hypersurface
\begin{equation}
M=\set{(x^0,x)}{x^0=u(x),\,x\in\mc S_0},
\end{equation}
then the induced metric has the form
\begin{equation}
g_{ij}=e^{2\psi}\{-u_iu_j+\s_{ij}\}
\end{equation}
where $\s_{ij}$ is evaluated at $(u,x)$, and its inverse $(g^{ij})=(g_{ij})^{-1}$ can
be expressed as
\begin{equation}\lae{1.10}
g^{ij}=e^{-2\psi}\{\s^{ij}+\frac{u^i}{v}\frac{u^j}{v}\},
\end{equation}
where $(\s^{ij})=(\s_{ij})^{-1}$ and
\begin{equation}\lae{1.11}
\begin{aligned}
u^i&=\s^{ij}u_j\\
v^2&=1-\s^{ij}u_iu_j\equiv 1-\abs{Du}^2.
\end{aligned}
\end{equation}
Hence, $\graph u$ is space-like if and only if $\abs{Du}<1$.

The covariant form of a normal vector of a graph looks like
\begin{equation}
(\n_\al)=\pm v^{-1}e^{\psi}(1, -u_i).
\end{equation}
and the contravariant version is
\begin{equation}
(\n^\al)=\mp v^{-1}e^{-\psi}(1, u^i).
\end{equation}
Thus, we have
\br Let $M$ be space-like graph in a future oriented coordinate system. Then, the
contravariant future directed normal vector has the form
\begin{equation}
(\n^\al)=v^{-1}e^{-\psi}(1, u^i)
\end{equation}
and the past directed
\begin{equation}\lae{1.15}
(\n^\al)=-v^{-1}e^{-\psi}(1, u^i).
\end{equation}
\er

In the Gau{\ss} formula \re{1.3} we are free to choose the future or past directed
normal, but we stipulate that we always use the past directed normal for reasons
that we have explained in \ci{cg5}.

Look at the component $\al=0$ in \re{1.3} and obtain in view of \re{1.15}

\begin{equation}\lae{1.16}
e^{-\psi}v^{-1}h_{ij}=-u_{ij}-\cha 000\mspace{1mu}u_iu_j-\cha 0j0
\mspace{1mu}u_i-\cha 0i0\mspace{1mu}u_j-\cha ij0.
\end{equation}
Here, the covariant derivatives a taken with respect to the induced metric of
$M$, and
\begin{equation}
-\cha ij0=e^{-\psi}\bar h_{ij},
\end{equation}
where $(\bar h_{ij})$ is the second fundamental form of the hypersurfaces
$\{x^0=\const\}$.

An easy calculation shows
\begin{equation}
\bar h_{ij}e^{-\psi}=-\tfrac{1}{2}\dot\s_{ij} -\dot\psi\s_{ij},
\end{equation}
where the dot indicates differentiation with respect to $x^0$.

Next, let us analyze under which condition a space-like hypersurface $M$ can be
written as a graph over the Cauchy hypersurface $\mc S_0$.

We first need
\bd

Let $M$ be a  closed, space-like hypersurface in $N$. Then,
\bi[(i)]
\item
$M$ is said to be \tit{achronal}, if no two points in $M$ can be connected by a
future directed time-like curve.

\item
$M$ is said to \tit{separate} $N$, if $N\raise 0.28ex\hbox{$\scriptstyle\setminus$}
M$ is disconnected.
\ei
\ed

In \ci[Proposition 2.5]{cg5} we proved

\bpp\lap{1.5}
Let $N$ be connected and globally hyperbolic, $\mc S_0\nobreak\ 
\su\nobreak\   N$ a compact Cauchy hypersurface, and $M\su N$ a compact,
connected space-like hypersurface of class $C^m, m\ge 1$. Then, $M=\graph \fu
u{\mc S}0$ with
$u\in C^m(\mc S_0)$ iff $M$ is achronal.
\epp

Sometimes, we need a Riemannian reference metric, e.g. if we want to estimate
tensors. Since the Lorentzian metric can be expressed as
\begin{equation}
\bar g_{\al\bet}dx^\al dx^\bet=e^{2\psi}\{-{dx^0}^2+\s_{ij}dx^i dx^j\},
\end{equation}
we define a Riemannian reference metric $(\tilde g_{\al\bet})$ by
\begin{equation}
\tilde g_{\al\bet}dx^\al dx^\bet=e^{2\psi}\{{dx^0}^2+\s_{ij}dx^i dx^j\}
\end{equation}
and we abbreviate the corresponding norm of a vectorfield $\h$ by
\begin{equation}
\nnorm \h=(\tilde g_{\al\bet}\h^\al\h^\bet)^{1/2},
\end{equation}
with similar notations for higher order tensors.

\section{The evolution problem}\las{2}

Let $N$ be a globally hyperbolic Lorentzian manifold with a compact Cauchy
hypersurface $\mc S_0$. Consider the problem of finding a closed hypersurface
of prescribed mean curvature $H$ in $N$, or more precisely, let $\Om$ be a
connected open subset of $N$, $f\in C^{0,\al}(\bar \Om)$, then we look for a
hypersurface $M\su \Om$ such that
\begin{equation}
\fv HM=f(x)\qq \A \,x\in M,
\end{equation}
where $\fv HM$ means that  $H$ is evaluated at the vector $(\ka_i(x))$ the
components of which are the principal curvatures of $M$.

We assume that $\pa \Om$ consists of two achronal, compact, connected,
space-like hypersurfaces $M_1$ and $M_2$, where $M_1$ is supposed to lie in the
\tit{past} of $M_2$. The $M_i$ should act as barriers for $(H,f)$.

\bd
$M_2$ is an \tit{upper barrier} for $(H,f)$, if $M_2$ is of class $C^{2,\al}$ and
\begin{equation}
\fv H{M_2}\ge f,
\end{equation}
and $M_1$ ia a \tit{lower barrier} for $(H,f)$, if $M_1$ is of class $C^{2,\al}$
satisfying
\begin{equation}
\fv H{M_1}\le f.
\end{equation}
\ed

In \ci[Section 6]{cg83} we proved the following theorem

\bt\lat{2.2}
Let $M_1$ be a lower and $M_2$ be an upper barrier for $(H,f)$, $f\in
C^{0,\al}(\bar \Om)$. Then, the problem
\begin{equation}\lae{2.4}
\fv HM=f
\end{equation}
has a solution $M\su \bar \Om$ of class $C^{2,\al}$ that can be written as a graph
over the Cauchy hypersurface $\mc S_0$.
\et

The crucial point in the proof is an a priori estimate in the $C^1$-norm and for
this estimate only the boundedness of $f$ is needed, i.e. even for merely
bounded $f$ $H^{2,p}$ solutions exist.

We want to give a new proof of \rt{2.2} that is based on the evolution method, and
to make this method work, we have to assume temporarily slightly higher
degrees of regularity for the barriers and right-hand side, i.e. we assume the
barriers to be of class $C^{4,\al}$ and $f$ to be of class $C^{2,\al}$. We can achieve
these assumptions by approximation without sacrificing the barrier conditions,
cf. \ci[p. 179]{cg97}.

To solve \re{2.4} we look at the evolution problem
\begin{equation}\lae{2.5}
\begin{aligned}
\dot x&=(H-f)\n,\\
x(0)&=x_0,
\end{aligned}
\end{equation}
where $x_0$ is an embedding of an initial hypersurface $M_0$, for which we
choose $M_0=M_2$,  $H$ is the mean curvature of the flow hypersurfaces
$M(t)$ with respect to the past directed normal $\n$, and $x(t)$ is an embedding
of $M(t)$.

In \ci{cg5} we have considered problems of the form \re{2.5} for general
curvature operators in a pseudo-riemannian setting, so that the present
situation can be retrieved as a special case of the general results in \ci[Section
3]{cg5}.

The evolution exists on a maximal time interval $[0,T^*)$, $0<T^*\le \un$, cf.
\ci[Section 2]{cg96}, where we apologize for the ambiguity of also calling the
evolution parameter \tit{time}.

Next, we want to show how the metric, the second fundamental form, and the
normal vector of the hypersurfaces $M(t)$ evolve. All time derivatives are
\tit{total} derivatives. We refer to \ci{cg5} for more general results and to
\ci[Section 3]{cg96}, where proofs are given in a Riemannian setting, but these
proofs are also valid in a Lorentzian environment.

\bl
The metric, the normal vector, and the second fundamental form of $M(t)$
satisfy the evolution equations
\begin{equation}
\dot g_{ij}=2(H- f)h_{ij},
\end{equation}
\begin{equation}\lae{2.7}
\dot \n=\nabla_M(H- f)=g^{ij}(H- f)_i x_j,
\end{equation}
and
\begin{equation}\lae{2.8}
\dot h_i^j=(H- f)_i^j- (H- f) h_i^k h_k^j-(H- f) \riema
\al\bet\ga\de\n^\al x_i^\bet \n^\ga x_k^\de g^{kj}
\end{equation}
\begin{equation}
\dot h_{ij}=(H- f)_{ij}+ (H- f) h_i^k h_{kj}-(H- f) \riema
\al\bet\ga\de\n^\al x_i^\bet \n^\ga x_j^\de.
\end{equation}
\el

\bl[Evolution of $(H- f)$]
The term $(H- f)$ evolves according to the equation
\begin{equation}\lae{2.10}
\begin{aligned}
{(H- f)}^\prime- \D (H- f)=&\msp[0]- \norm A^2(H- f)-f_\al\n^\al (H- f)\\
&-\bar R_{\al\bet}\n^\al\n^\bet (H- f),
\end{aligned}
\end{equation}
where
\begin{equation}
(H- f)^{\prime}=\frac{d}{dt}(H- f)
\end{equation}
and
\begin{equation}
\norm A^2=h_{ij}h^{ij}.
\end{equation}
\el

From \re{2.8} we deduce with the help of the Ricci identities a parabolic equation
for the second fundamental form
\bl\lal{2.7}
The mixed tensor $h_i^j$ satisfies the parabolic equation
\begin{equation}\raisetag{-58pt}\lae{2.13}
\begin{aligned}
\dot h_i^j-\D h_i^j&=-\norm A^2h_i^j+f h_i^kh_k^j
- f_{\al\bet} x_i^\al x_k^\bet g^{kj}- f_\al\n^\al h_i^j\\
&\q\,+2\riema \al\bet\ga\de x_m^\al x_i ^\bet x_k^\ga
x_r^\de h^{km} g^{rj}\\
&\q\,-g^{kl}\riema \al\bet\ga\de x_m^\al x_k ^\bet x_r^\ga x_l^\de
h_i^m g^{rj}- g^{kl}\riema \al\bet\ga\de x_m^\al x_k ^\bet x_i^\ga x_l^\de h^{mj} \\
&\q\,-\bar R_{\al\bet}\n^\al\n^\bet h_i^j+f\riema \al\bet\ga\de\n^\al x_i^\bet\n^\ga x_m^\de
g^{mj}\\ 
&\q\,+ g^{kl}\bar R_{\al\bet\ga\de;\e}\{\n^\al x_k^\bet x_l^\ga x_i^\de
x_m^\e g^{mj}+\n^\al x_i^\bet x_k^\ga x_m^\de x_l^\e g^{mj}\}.
\end{aligned}
\end{equation}
\el

\br\lar{2.6}
In view of the maximum principle, we immediately deduce from \re{2.10} that the
term $(H-f)$ has a sign during the evolution if it has one at the beginning.
Thus, we have
\begin{equation}\lae{2.14}
H\ge f.
\end{equation}
\er

\section{Lower order estimates}\las 3

We recall our  assumption that the ambient space is
globally hyperbolic with a compact Cauchy hypersurface $\so$. The barriers
$M_i$ are then graphs over $\so, M_i=\graph u_i$, because they are achronal, cf.
\rp{1.5}, and we have
\begin{equation}\lae{3.1}
u_1\le u_2,
\end{equation}
for $M_1$ should lie in the past of $M_2$, and the enclosed domain is supposed to
be connected. Moreover, in view of the Harnack inequality, the strict inequality
is valid in \re{3.1} unless the barriers coincide and are a solution to our
problem.

Let us look at the evolution equation \re{2.5} with initial hypersurface $M_0$
equal to $M_2$. Then, because of the short-time existence, the
evolution will exist on a maximal time interval
$I=[0,T^*), \,T^*\le
\un$, as long as the evolving hypersurfaces are space-like and
smooth.

Furthermore, since the initial hypersurface is a graph over $\so$, we can write
\begin{equation}
M(t)=\graph\fu{u(t)}S0\q \A\,t\in I,
\end{equation}
where $u$ is defined in the cylinder $Q_{T^*}=I\times \so$. We then deduce from
\re{2.5}, looking at the component $\al=0$, that $u$ satisfies a parabolic equation
of the form
\begin{equation}\lae{3.3}
\dot u=-e^{-\psi}v^{-1}(H-f),
\end{equation}
where we  use the notations in \rs{1}, and where we emphasize that the time
derivative is a total derivative, i.e.
\begin{equation}\lae{3.4}
\dot u=\pde ut+u_i\dot x^i.
\end{equation}

Since the past directed normal can be expressed as
\begin{equation}
(\n^\al)=-e^{-\psi}v^{-1}(1,u^i),
\end{equation}
we conclude from \re{2.5}, \re{3.3}, and \re{3.4}
\begin{equation}
\lae{3.6}
\pde ut=-e^{-\psi}v(H- f).
\end{equation}
Thus, $\pde ut$ is non-positive in view of \rr{2.6}.

Next, let us state our first a priori estimate

\bl\lal{3.1}
During the evolution the flow hypersurfaces stay in $\bar \Om$.
\el

This is an immediate consequence of the Harnack inequality, cf. \ci[Lemma
5.1]{cg96} for details.

As a consequence of \rl{3.1} we obtain
\begin{equation}
\inf_{\so} u_1\le u\le \sup_\so u_2\q \A\,t\in I.
\end{equation}

We are now able to derive the $C^1$-estimates, i.e. we shall show that the
hypersurfaces remain uniformly space-like, or equivalently, that the term
\begin{equation}
\tilde v=v^{-1}=\frac{1}{\sql}
\end{equation}
is uniformly bounded.

Let us first derive an evolution equation for $\tilde v$.

\bl[Evolution of $\tilde v$]\lal{3.2}
Consider the flow \re{2.5} in the distinguished coordinate system associated
with $\so$. Then, $\tilde v$ satisfies the evolution equation
\begin{equation}\lae{3.20}
\begin{aligned}
\dot{\tilde v}-\D\tilde v=&-\norm A^2\tilde v
-f\h_{\al\bet}\n^\al\n^\bet-f_\bet x_i^\bet  \h_\al x_k^\al g^{ik}\\
&-2h^{ij} x_i^\al x_j^\bet \h_{\al\bet}-g^{ij}\h_{\al\bet\ga}x_i^\bet
x_j^\ga\n^\al\\
&-\bar R_{\al\bet}\n^\al x_k^\bet\h_\ga x_l^\ga g^{kl},
\end{aligned}
\end{equation}
where $\h$ is the covariant vector field $(\h_\al)=e^{\psi}(-1,0,\dotsc,0)$.
\el

\bp
We have $\tilde v=\spd \h\n$. Let $(\x^i)$ be local coordinates for $M(t)$.
Differentiating $\tilde v$ covariantly we deduce
\begin{equation}\lae{3.21}
\tilde v_i=\h_{\al\bet}x_i^\bet\n^\al+\h_\al\n_i^\al,
\end{equation}
\begin{equation}\lae{3.22}
\begin{aligned}
\tilde v_{ij}= &\msp[5]\h_{\al\bet\ga}x_i^\bet x_j^\ga\n^\al+\h_{\al\bet}x_{ij}^\bet\n^\al\\
&+\h_{\al\bet}x_i^\bet\n_j^\al+\h_{\al\bet}x_j^\bet\n_i^\al+\h_\al\n_{ij}^\al
\end{aligned}
\end{equation}

The time derivative of $\tilde v$ can be expressed as
\begin{equation}\lae{3.23}
\begin{aligned}
\dot{\tilde v}&=\h_{\al\bet}\msp\dot x^\bet\n^\al+\h_\al\dot\n^\al\\
&=\h_{\al\bet}\n^\al\n^\bet(H-f)+(H-f)^k x_k^\al\h_\al\\
&=\h_{\al\bet}\n^\al\n^\bet(H-f)+H^k x_k^\al\h_\al-{ f}_\bet x_i^\bet
x_k^\al g^{ik}\h_\al,
\end{aligned}
\end{equation}
where we have used \re{2.7}.

Substituting \re{3.22} and \re{3.23} in \re{3.20}, and simplifying the resulting
equation with the help of the Weingarten and Codazzi equations, we arrive at the
desired conclusion.
\ep

\bl\lal{3.3}
There is a constant $c=c(\Om)$ such that for any positive function $0<\e=\e(x)$ on
$\so$ and any hypersurface $M(t)$ of the flow we have
\begin{align}
\nnorm \n&\le c\tilde v,\\\lae{3.14}
g^{ij}&\le c\tilde v^2\s^{ij},\\
\intertext{and}\lae{3.15}
\abs{h^{ij}\h_{\al\bet}x_i^\al x_j^\bet}&\le \frac{\e}{2}\norm A^2\tilde
v+\frac{c}{2\e}\tilde v^3
\end{align}
where $(\h_\al)$ is the vector field in \rl{3.2}.
\el

\bp
The first two estimates can be immediately verified. To prove \re{3.15} we
choose local coordinates $(\x^i)$ such that
\begin{equation}
h_{ij}=\ka_i\de_{ij},\qq g_{ij}=\de_{ij}
\end{equation}
and deduce
\begin{equation}
\begin{aligned}
\abs{h^{ij}\h_{\al\bet}x_i^\al x_j^\bet}&\le \sum_i\abs{\ka_i}\abs{\h_{\al\bet}x_i^\al
x_i^\bet}\\
&\le \frac{\e}{2}\norm A^2\tilde v+\frac{1}{2\e}\tilde v^{-1}\sum_i\abs{\h_{\al\bet}
x_i^\al x_i^\bet}^2,
\end{aligned}
\end{equation}
and
\begin{equation}
\sum_i\abs{\h_{\al\bet}
x_i^\al x_i^\bet}^2\le g^{ik}\h_{\al\bet}x_i^\al x_j^\bet \msp[3]g^{jl} \h_{\ga\de} x_k^\ga
x_l^\de.
\end{equation}
Hence, the result in view of \re{3.14}.
\ep

Combining the preceding lemmata we infer

\bl\lal{3.4}
There is a constant $c=c(\Om)$ such that for any positive function $\e=\e(x)$ on
$\so$ the term $\tilde v$ satisfies a parabolic inequality of the form
\begin{equation}
\dot{\tilde v}-\D\tilde v\le -(1-\e)\norm A^2\tilde v+c[\abs f+\nnorm{Df}]\tilde
v^2+c[1+\e^{-1}]\tilde v^3.
\end{equation}
\el

We note that the statement \tit{$c$ depends on $\Om$} also implies that $c$
depends on geometric quantities of the ambient space restricted to $\Om$.

We further need the following two lemmata

\bl\lal{3.5}
Let $M(t)=\graph u(t)$ be the flow hypersurfaces, then we have
\begin{equation}
\dot u-\D u=e^{-\psi}v^{-1}f-e^{-\psi}g^{ij}\bar h_{ij}+\cha 000\norm{Du}^2+2\cha
0i0 u^i,
\end{equation}
where the time derivative is a total derivative.
\el

\bp
We use the relation
\begin{equation}
\dot u=-e^{-\psi}v^{-1}(H-f)
\end{equation}
together with \re{1.16}.
\ep

\bl\lal{3.6}
Let $M\su \bar \Om$ be a graph over $\so$, $M=\graph u$, then
\begin{equation}
\abs{\tilde v_i u^i}\le c\tilde v^3+\norm A e^\psi\norm {Du}^2,
\end{equation}
where $c=c(\Om)$.
\el

\bp
First, we use that
\begin{equation}
\tilde v^2=1+e^{2\psi}\norm{Du}^2,
\end{equation}
and thus,
\begin{equation}
2\tilde v\tilde v_i=2\psi_\al x_i^\al e^{2\psi} \norm{Du}^2+2e^{2\psi} u_{ij}u^j,
\end{equation}
from which we infer
\begin{equation}
\abs{\tilde v_iu^i}\le c\tilde v^3+\tilde v^{-1}e^{2\psi}\abs{u_{ij}u^i u^j},
\end{equation}
which gives the result because of \re{1.16}.
\ep

We are now ready to prove the uniform boundedness of $\tilde v$.

\bpp
During the evolution the term $\tilde v$ remains uniformly bounded
\begin{equation}
\tilde v\le c=c(\Om,\abs f,\nnorm{Df}).
\end{equation}
\epp

\bp
Let $\m,\lam$ be positive constants, where $\m$ is supposed to be small and $\lam$
large, and define
\begin{equation}\lae{3.27}
\f=e^{\m e^{\lam u}},
\end{equation}
where we assume without loss of generality that $1\le u$, otherwise replace in
\re{3.27} $u$ by $(u+c)$, $c$ large enough.

We shall show that
\begin{equation}
w=\tilde v \f
\end{equation}
is uniformly bounded if $\m,\lam$ are chosen appropriately.

In view of \rl{3.3} and \rl{3.5} we have
\begin{equation}
\dot\f-\D\f\le c\m\lam e^{\lam u}[\tilde v\abs f +\tilde v^2] \f-\m\lam^2 e^{\lam u} [1+\m
e^{\lam u}]\norm{Du}^2\f,
\end{equation}
from which we further deduce taking \rl{3.4} and \rl{3.6} into account
\begin{equation}
\begin{aligned}
\dot w-\D w&\le -(1-\e) \norm A^2\tilde v\f +c[\abs f+\nnorm{Df}]\tilde v^2\f\\
&\q\,+c[1+\e^{-1}]\tilde v^3\f-\m\lam^2 e^{\lam u} [1+\m e^{\lam u}] \tilde v
\norm{Du}^2\f\\
&\q\,+c[1+\abs f]\m\lam e^{\lam u}\tilde v^3\f+2\m\lam e^{\lam u} \norm A e^\psi
\norm{Du}^2\f.
\end{aligned}
\end{equation}

We estimate the last term on the right-hand side by

\begin{equation}
\begin{aligned}
2\m\lam e^{\lam u}\norm A e^\psi\norm{Du}^2\f&\le (1-\e)\norm A^2\tilde v\f\\
&\q\,+\frac{1}{1-\e}\m^2\lam^2e^{2\lam u}\tilde v^{-1}e^{2\psi}\norm{Du}^4\f,
\end{aligned}
\end{equation}
and conclude

\begin{equation}
\begin{aligned}
\dot w-\D w&\le  c[\abs f+\nnorm{Df}]\tilde v^2\f+ c[1+\abs f]\m\lam e^{\lam u} \tilde
v^3\f\\
 &\q\,+c[1+\e^{-1}]\tilde v^3\f +[\frac{1}{1-\e}-1]\m^2\lam^2 e^{2\lam
u}\norm{Du}^2\tilde v\f\\
&\q\,-\m\lam^2 e^{\lam u}\norm{Du}^2\tilde v\f,
\end{aligned}
\end{equation}
where we have used that
\begin{equation}
e^{2\psi}\norm{Du}^2\le \tilde v^2.
\end{equation}

Setting $\e=e^{-\lam u}$, we then obtain

\begin{equation}\lae{3.34}
\begin{aligned}
\dot w-\D w&\le c[\abs f+\nnorm{Df}]\tilde v^2\f+c e^{\lam u} \tilde v^3\f\\
&\q\,+c[1+\abs f]\m\lam e^{\lam u}\tilde v^3\f\\
&\q\,+[\frac{\m}{1-\e}-1]\m\lam^2 e^{\lam u}\norm{Du}^2\tilde v\f.
\end{aligned}
\end{equation}

Now, we choose $\m=\frac{1}{2}$ and $\lam_0$ so large that
\begin{equation}
\frac{\m}{1-e^{-\lam u}}\le \frac{3}{4}\qq\A\,\lam\ge \lam_0,
\end{equation}
and infer that the last term on the right-hand side of \re{3.34} is less than
\begin{equation}
-\frac{1}{8}\lam^2e^{\lam u}\norm{Du}^2\tilde v\f
\end{equation}
which in turn can be estimated from above by
\begin{equation}
-c\lam^2e^{\lam u}\tilde v^3\f
\end{equation}
at points where $\tilde v\ge 2$.

Thus, we conclude that for
\begin{equation}
\lam\ge \max (\lam_0, 4[1+\abs f_{_\Om}])
\end{equation}
the parabolic maximum principle, applied to $w$, yields
\begin{equation}
w\le \const (\abs{w(0)}_{_\so},\lam_0,\abs f, \nnorm{Df},\Om).
\end{equation}
\ep

\section{$C^2$-estimates}\las{4}

Since the mean curvature operator is a quasilinear operator, the uniform
$C^1$-estimates we have established in the last section also yield uniform
$C^2$-estimates during the evolution, but nevertheless, we would like to give an
independent proof of the $C^2$-estimates.

\bl\lal{4.1}
During the evolution the principal curvatures of the evolution
hypersurfaces $M(t)$ are uniformly bounded.
\el

\bp
As already mentioned in \rr{2.6}, we know that $f\le H$, thus, it is sufficient to
estimate the principal curvatures from above. 

Let $\f$  be defined  by
\begin{equation}\lae{4.1}
\f=\sup\set{{h_{ij}\h^i\h^j}}{{\norm\h=1}}.
\end{equation}
We claim that
$\f$ is uniformly bounded.

Let $0<T<T^*$, and $x_0=x_0(t_0)$, with $ 0<t_0\le T$, be a point in $M(t_0)$ such
that

\begin{equation}
\sup_{M_0}\f<\sup\set {\sup_{M(t)} \f}{0<t\le T}=\f(x_0).
\end{equation}

We then introduce a Riemannian normal coordinate system $(\x^i)$ at $x_0\in
M(t_0)$ such that at $x_0=x(t_0,\x_0)$ we have
\begin{equation}
g_{ij}=\de_{ij}\q \tup{and}\q \f=h_n^n.
\end{equation}

Let $\tilde \h=(\tilde \h^i)$ be the contravariant vector field defined by
\begin{equation}
\tilde \h=(0,\dotsc,0,1),
\end{equation}
and set
\begin{equation}
\tilde \f=\frac{h_{ij}\tilde \h^i\tilde \h^j}{g_{ij}\tilde \h^i\tilde \h^j}\raise 2pt
\hbox{.}
\end{equation}

$\tilde \f$ is well defined in neighbourhood of $(t_0,\x_0)$, and $\tilde \f$
assumes its maximum at $(t_0,\x_0)$. Moreover, at $(t_0,\x_0)$ we have
\begin{equation}
\dot{\tilde \f}=\dot h_n^n,
\end{equation}
and the spatial derivatives do also coincide; in short, at $(t_0,\x_0)$ $\tilde \f$
satisfies the same differential equation \re{2.13} as $h_n^n$. For the sake of
greater clarity, let us therefore treat $h_n^n$ like a scalar and pretend that
$\f=h_n^n$. 

At $(t_0,\x_0)$ we have $\dot\f\ge 0$, and, in view of the maximum principle, we
deduce from \rl{2.7}
\begin{equation}
0\le -\norm A^2h_n^n+f\abs{h_n^n}^2+c[\abs
f+\nnorm{Df}+\nnorm{D^2f}][1+\abs{h_n^n}].
\end{equation}

Thus, $\f$ is uniformly bounded.
\ep

\section{Convergence to a stationary solution}\las 5

We are now ready to give a new proof of  \rt{2.2}.  Let us look at the scalar
version of the flow as in \re{3.6}
\begin{equation}\lae{5.1}
\pde ut=-e^{-\psi}v(H- f).
\end{equation}
This is  a scalar parabolic differential equation defined on the cylinder
\begin{equation}
Q_{T^*}=[0,T^*)\times \so
\end{equation}
with initial value $u(0)=u_2\in C^{4,\al}(\so)$. In view of the a priori estimates,
which we have established in the preceding sections, we know that
\begin{equation}
{\abs u}_\low{2,0,\so}\le c
\end{equation}
and
\begin{equation}
H\,\tup{is uniformly elliptic in}\,u
\end{equation}
independent of $t$. Thus, we can apply
the known regularity results, see e.g. \ci[Chapter 5.5]{nk}, where even more
general operators are considered,  to conclude that uniform
$C^{2,\al}$-estimates are valid, leading further to uniform $C^{4,\al}$-estimates
due to the regularity results for linear operators.

Therefore, the maximal time interval is unbounded, i.e. $T^*=\un$.

Now, integrating \re{5.1} with respect to $t$, and observing that the right-hand
side is non-positive, yields
\begin{equation}
u(0,x)-u(t,x)=\int_0^te^{-\psi}v(H- f)\ge c\int_0^t(H- f),
\end{equation}
i.e.,
\begin{equation}
\int_0^\un \abs{H- f}<\un\qq\A\msp x\in \so
\end{equation}
Hence, for any $x\in\so$ there is a sequence $t_k\rightarrow \un$ such that
$(H- f)\rightarrow 0$.

On the other hand, $u(\cdot,x)$ is monotone decreasing and therefore
\begin{equation}
\lim_{t\rightarrow \un}u(t,x)=\tilde u(x)
\end{equation}
exists and is of class $C^{4,\al}(\so)$ in view of the a priori estimates. We, finally,
conclude that $\tilde u$ is a stationary solution of our problem, and that
\begin{equation}
\lim_{t\rightarrow \un}(H- f)=0.
\end{equation}

To prove existence under the weaker assumptions of \rt{2.2}, we use
approximation and  the a priori estimate in \ci[Theorem 4.1]{cg83}.


\end{document}